\documentclass[12pt]{article}

\newtheorem{theorem}{Theorem}
\def \proof {\par \noindent {\bf Proof.}\hskip 5pt}

\begin{document}

\title{On empty pentagons and hexagons in planar point sets\footnote{Work on
this paper was supported by project 1M0545 of
The Ministry of Education of the Czech Republic.}}
\author{Pavel Valtr\footnote{Department of Applied Mathematics and Institute for Theoretical Computer Science (ITI),
Charles University, Malostransk\'e n\'am.~25, 182 00 Praha 1, Czech Republic}
%Email:~{\tt roddick@csem.flinders.edu.au}\\[.1in]
%$^2$ Some other Department \\
%Some other University, \\
%Somewhere else, \\
%Email:~{\tt Anne.Other@cs.somewhere.edu}
}

\date{}

\maketitle

%Include either the full copyright statement:
%\toappear{Copyright \copyright 2006, Australian Computer Society, Inc.  This paper appeared at the Seventeenth Australasian Database Conference (ADC2006), Hobart, Australia.  Conferences in Research and Practice in Information Technology (CRPIT), Vol. 49. Gillian Dobbie and James Bailey, Eds. Reproduction for academic, not-for profit purposes permitted provided this text is included. }

%or more easily (and recommended) use the alternative:

%  For Government work where the alternative copyright statement will be required, please replace the \toappearstandard clause by (uncommented)
%\newcommand\copyrightholder{Commonwealth of Australia} - or whoever
%\toappeargovwork

\begin{abstract}
We give improved lower bounds on the minimum number of $k$-holes (empty convex $k$-gons)
in a set of $n$ points in general position in the plane, for $k=5,6$.
 
\end{abstract}
\vspace{.1in}


\section{Introduction}

We say that a set $P$ of points in the plane is in {\it general position} if
it contains no three points on a line.

Let $P$ be a set of $n$ points in general position in the plane.
A {\em k-hole of $P$} (sometimes also called empty convex $k$-gon or convex $k$-hole)
is a set of vertices of a convex $k$-gon with vertices in $P$ containing no other points of $P$.

Let $X_k (n)$ be the minimum number of $k$-holes in a set of $n$
points in general position in the plane.
Horton \cite{ho} proved that $X_k(n)=0$ for any $k\ge7$
and for any positive integer $n$.
The following bounds on $X_k (n), k=3,4,5,6$, are known (the letter
$H$ denotes the number of vertices of the convex hull of the point set):

$$n^2-\frac{43}{9}n+H+\frac{11}{9}\le X_3(n) \le 1.6195...n^2+o(n^2),$$

$$\frac{n^2}{2}-\frac{55}{18}n+H-\frac{23}{9}\le  X_4(n) \le 1.9396...n^2+o(n^2),$$
$$\frac{2}{9}n-O(1) \le X_5(n) \le 1.0206... n^2+o(n^2),$$

$$ \frac{n}{463} - 1  \le X_6(n) \le 0.2005...n^2+o(n^2).$$
The upper bounds were shown in \cite{bv},
improving previous bounds of \cite{km,bf,v2,d}. The lower bounds for $k=3,4,5$
can be found in an updated version of the conference paper \cite{ga},
also improving lower bounds from several papers.
The lower bound on $X_6(n)$ follows from a result of V.~A.~Koshelev~\cite{ko}.
In this paper we give the following improved lower bounds:
%%%%%%%%%%%%%%%%%%%%%%%%%%%%%%%%%%
\begin{theorem}\label{t:main}
\begin{eqnarray*}
X_5(n)  \ge n/2 - O(1) ,\\
X_6(n) \ge  n/229-4.
\end{eqnarray*}
\end{theorem}

After finishing our research, we have learned that a group of researchers
including Oswin Aichholzer, Ruy Fabila-Monroy, Clemens Huemmer, and Birgit Vogtenhuber
has very recently obtained a better bound $X_5(n)  \ge 3n/4 - o(n)$.
Their result is not written yet. Their method does not seem to achieve
our bound on $X_6(n)$ but it also gives slight
improvements on the lower bounds on $X_3(n)$ and $X_4(n)$ mentioned above.

\section{Proofs}
To prove the first inequality in Theorem~\ref{t:main}, it suffices to prove that
if $P$ is a set of $n>20$ points in general position in the plane
then $P$ contains a subset $P'$ of eight points such that $P'$ and $P-P'$ can be
separated by a line and at least four $5$-holes of $P$ intersect $P'$.
Indeed, if this is true then we can repeatedly remove eight points of $P'$.
Each removal decreases the number of points by $8$ and the number of $5$-holes by at least $4$.
Doing this as long as at least $21$ points remain, we obtain the first inequality in Theorem~\ref{t:main}.

Let $P$ be a set of $n>20$ points in general position in the plane.
For two points $x,y$ of $P$, we denote by $L(xy)$ the open halfplane to the left of the line $xy$ (oriented from $x$ to $y$). The
complementary open halfplane is denoted by $R(xy)$. If $L(xy)$ contains exactly $k$ points of $P$, then we say that the oriented
segment $xy$ is a $k$-edge of $P$.

Take a vertex $a$ of the convex hull of $P$. Order the other points radially around $a$ starting from 
the point on the convex hull clockwise from $a$. Let $a'$ be the 12-th point in this order. Then $aa'$ is an
$11$-edge. Since $X_5(10)>0$ \cite{harb},  $L(aa')$ contains a $5$-hole, $D$, of $P$. In the rest of the proof,
$D$ is fixed but $aa'$ may later denote other $11$-edges.

The key part of the proof is to find an $11$-edge $bb'$ such that $b$ is a vertex of $D$ and the other four vertices of $D$ lie in
$L(bb')$. To do it, we clockwise rotate a line $l$ starting from $l=aa'$ as follows. Initially
we start to rotate $l$ at the midpoint of the segment $aa'$. During the rotation, the center of
rotation may change at any moment but the rotated line $l$ cannot go over any point of $P$. We rotate as long as it is
possible, until we reach a position $l=bb'$, where $b,b' \in P$, the point $b$ was originally to the left of $l$ and $b'$
was originally to the right of $l$.
Thus, $b \in L(aa')\cup\{a'\}$ and $b' \in R(aa')\cup\{a\}$. There are no points of $P$ in the open wedges $R(aa') \cap L(bb')$ and
$L(aa') \cap R(bb')$. The edge $bb'$ is an $11$-edge of $P$. We distinguish three cases:

Case 1: The segments $aa'$ and $bb'$ internally cross, thus $a,a',b,b'$ are pairwise different.

Case 2: $b'=a$.

Case 3: $b=a'$.

Since $D$ lies in $L(aa')$, it also lies in $L(bb')\cup\{b\}$. The point $b$ may be a vertex of $D$ in Cases 1 and 2.
All other vertices of $D$ lie in $L(bb')$. If $b$ is not a vertex of $D$, then we rename the points $b$ and $b'$ by $a$ and $a'$,
respectively, and rotate a line $l$ in the same way as above from the position $l=aa'$. We reach some new position $l=bb'$.
Repeat this process until the point $b$ coincides with one of the vertices of $D$. (Note that the line $l$ cannot
rotate outside of $D$ forever, because $n>20$.) Then we are in Case 1 or in Case 2, and the
other four vertices of $D$ lie in $L(aa') \cap L(bb')$. In Case 1 or 2, we consider the $12$-point set $Q:= (P \cap
L(bb')) \cup \{b'\}$. Since $X_5(12)\ge 3$ \cite{dehn}, the set $Q$
contains at least three $5$-holes of $P$. Together with $D$, these are at least four $5$-holes of $P$ with
vertices in the $13$-point set $Q \cup \{b\} = P \cap closure(L(bb'))$. None of these $5$-holes contains both $b$ and $b'$.
Therefore, we can take $P'$ as the set of eight points of $L(bb')$ with largest distances to the line $bb'$.
This finishes the proof of the first inequality in Theorem~\ref{t:main}.
\medskip

We remark without proof that a slightly better bound $(1/2+c)n - const$ with $c>0$ can be obtained by using the fact that 
any sufficiently large set $P$ contains linearly many disjoint $6$-holes.

The above proof can be generalized to give the more general theorem below.
The theorem below together with $X_6(463)>0$ (proved by V.~A.~Koshelev~\cite{ko})
gives the second inequality in Theorem~\ref{t:main}.

\begin{theorem}\label{t:general}
Suppose that $X_k(s-1)\ge 1$ and $X_k(s)\ge t$ for some positive integers $k,s,t$.
Then $X_k(n)\ge \frac{t+1}{s-k+1}(n-(2s-2))$ for $n\ge 2s-2$.
\end{theorem}

%\bigskip

%{\em Proof.} 
\proof If $P$ is a set of $n>2s-2$ points
then $P$ contains an ($s-1$)-edge $aa'$. Let $D$ be a $k$-hole of $P$ contained in $L(aa')$.
Analogously as in the previous proof, we find two $(s-1)$-edges $aa'$ and $bb'$ such that
$b$ is a vertex of $D$ and $D$ lies in $L(aa')$ and also in $L(bb')\cup\{b\}$. 
In Case 1 or 2, we consider the $s$-point set $Q:= (P \cap
L(bb')) \cup \{b'\}$. Since $X_k(s) \ge t$,
the set $Q$ contains at least $t\ k$-holes of $P$. Together with $D$, these are at least $t+1\ k$-holes of $P$ with
vertices in the $s+1$-point set $Q \cup \{b\} = P \cap closure(L(bb'))$. None of these $k$-holes contains both $b$ and $b'$.
Therefore, if we take $P'$ as the set of $s-k+1$ points of $L(bb')$ with largest distances to the line $bb'$ then removing the $s-k+1$ points of $P'$ from $P$ decreases the number of $k$-holes
by at least $t+1$.
Theorem~\ref{t:general} follows.

%%%%%%%%%%%%%%%%%%%%%%%%%%%%%%

%%%%%%%%%%%%%%%%%%%%%

\end{document}